\documentclass[10pt,reqno]{amsart} 

\usepackage{amssymb,latexsym}
\usepackage{cite} 

\usepackage[height=190mm,width=130mm]{geometry} 

\usepackage{amsfonts}
\usepackage{amsmath}
\usepackage{amscd}
\usepackage{amssymb}
\usepackage{amsthm}
\usepackage{bm}
\usepackage{enumerate}
\usepackage{enumitem}
\usepackage{hyperref}
\usepackage{latexsym}
\usepackage{mathtools}

\newcommand{\tpi}{2\pi i}

\newcommand{\ZZ}{\mathcal Z}

\newcommand{\C}{\mathbb{C}}

\newcommand{\I}{\mathcal{I}}

\newcommand{\Ip}{\I_{+}}

\newcommand{\N}{\mathbb{N}}

\newcommand{\wt}{\textup{wt}}
\newcommand{\Z}{\mathbb{Z}}

\newcommand{\R}{\mathbb{R}}
\newcommand{\HH}{\mathbb{H}}

\usepackage{hyperref}

\theoremstyle{plain}
\newtheorem{theorem}{Theorem}[section]
\newtheorem{lemma}[theorem]{Lemma}

\newtheorem{proposition}[theorem]{Proposition}

\theoremstyle{definition}
\newtheorem{definition}[theorem]{Definition}

\theoremstyle{remark}
\newtheorem{remark}[theorem]{Remark}

\DeclareMathOperator{\Ker}{Ker}
\DeclareMathOperator{\Ima}{Im}

\numberwithin{equation}{section} 

\newcommand{\g}{\ensuremath{\Gamma}}

\newcommand{\ps}{{\raise 1pt\hbox{\tiny (}}}

\newcommand{\pss}{{\raise 1pt\hbox{\tiny [}}}
\newcommand{\pdd}{{\raise 1pt\hbox{\tiny ]}}}
\newcommand{\pd}{{\raise 1pt\hbox{\tiny )}}}

\newcommand{\bs}{{\raise 1pt\hbox{\tiny [}}}
\newcommand{\bd}{{\raise 1pt\hbox{\tiny ]}}}

\def\cross{\mathinner{\mathrel{\raise0.8pt\hbox{$\scriptstyle>$}}
                 \joinrel\mathrel\triangleleft}}

\usepackage{amsmath,amsthm}
\usepackage{amsfonts}
\usepackage{amssymb}
\usepackage{longtable}
\usepackage[matrix,arrow,curve]{xy}

\def\K{\mathcal{K}}

\usepackage{stackrel}

\newcommand{\be}{\begin{equation}}
\newcommand{\ee}{\end{equation}}

\newcommand{\nc}{\newcommand}
\nc{\cali}{\mathcal}
\nc{\on}{\operatorname}
\nc{\Wick}{{\mb :}}

\nc{\ddz}{\frac{\partial}{\partial z}}
\nc{\ch}{\mbox{ch}}
\nc{\Oo}{{\cali O}}
\nc{\cond}{|\,}
\nc{\bib}{\bibitem}
\nc{\pone}{\Pro^1}
\nc{\pa}{\partial}
\nc{\arr}{\rightarrow}
\nc{\larr}{\longrightarrow}
\nc{\ket}{\rangle}
\nc{\bra}{\langle}
\nc{\gam}{\bar{\gamma}}
\nc{\q}{\widetilde{Q}}
\nc{\ep}{\epsilon}
\nc{\su}{\widehat{{\mf s}{\mf l}}_2}
\nc{\sw}{{\mf s}{\mf l}}
\nc{\h}{{\mf h}}
\nc{\n}{{\mf n}}
\nc{\ab}{\mf{a}}
\nc{\is}{{\mb i}}
\nc{\js}{{\mb j}}
\nc{\bi}{\bibitem}
\nc{\He}{{\cali H}}
\nc{\inv}{^{-1}}
\nc{\ol}{\overline}
\nc{\wh}{\widehat}
\nc{\dst}{\displaystyle}

\nc{\delt}{\partial_t}
\nc{\ddt}{\frac{\partial}{\partial t}}
\nc{\delx}{\partial_x}
\nc{\mb}{\mathbf}
\nc{\mf}{\mathfrak}

\nc{\mbb}{\mathbb}
\nc{\Ctt}{\C((t))}
\nc{\Ct}{\C[t,t\inv]}

\nc{\ghat}{\wh{\g}}

\nc{\un}{\underline}
\nc{\mc}{\mathcal}
\nc{\BB}{{\mc B}}
\nc{\bb}{{\mf b}}
\nc{\kk}{{\mf k}}
\nc{\frob}{\times}
\nc{\sm}{\setminus}
\nc{\Pp}{{\mathbb P}^1}
\nc{\Aa}{{\mc A}}

\nc{\AutO}{\on{Aut}\Oo}
\nc{\AUTO}{\un{\on{Aut}}\Oo}
\nc{\AUTK}{\un{\on{Aut}}\K}
\nc{\Heout}{\He_{\out}}
\nc{\Hetil}{{\widetilde\He}}
\nc{\wb}{\overline}

\nc{\Res}{\on{Res}}
\nc{\pitil}{\Pi}
\nc{\Ctil}{\wt{C}}
\nc{\auto}{\on{Aut} \Oo}
\nc{\phitil}{\wt{\phi}}
\nc{\gz}{\g_{\vec z}}
\nc{\tensorM}{\bigotimes_{i=1}^N{\mathbb M}_i}
\nc{\tensorW}{\bigotimes_{i=1}^N W_{\nu_i,k}}
\nc{\out}{\on{out}}

\nc{\m}{{\mathfrak m}}

\nc{\gx}{\g^0_{\vec x}}
\nc{\hx}{\He^0_{\vec x}}
\nc{\tensorpi}{\pi_{\nu_1,\ldots,\nu_N}^\kappa}
\nc{\Phizw}{\Phi_{\vec w}({\vec z})}
\nc{\Pro}{{\mathbb P}}

\nc{\De}{\Delta}

\nc{\us}{\underset}

\nc{\Ll}{\mc L}
\nc{\dR}{\on{dR}}

\nc{\T}{{\mc T}}

\nc{\Xn}{\overset{\circ}X{}^n} \nc{\Dn}{\overset{\circ}D{}^n}
\nc{\Dxn}{\overset{\circ}D{}^n_x} \nc{\varphitil}{\wt{\varphi}}

\nc{\lf}{{\mf l}}
\nc{\GL}{{}^L G}
\nc{\Vir}{\on{Vir}}

\begin{document}
\title[Cohomology of multipoint connections on complex curves]
{Cohomology of multipoint connections on complex curves}
\author{A. Zuevsky}
\address{Institute of Mathematics \\ Czech Academy of Sciences\\ Zitna 25, Prague \\ Czech Republic}
\email{zuevsky@yahoo.com}

\begin{abstract}
Given complex functions on complex curves satisfying recursion relations
with respect to the number of marked points at which they are evaluated,
we construct a cohomology theory governing such recursions, expressed in
terms of a generalization of holomorphic connections which we call
multipoint connections. We identify precisely the sense in which the
associated coboundary operators square to zero, as an integrability
condition on the recursion kernels rather than as an identity forcing the
correlation functions themselves to vanish, and prove, following the
general theory of systems of functional equations, that this condition
has non-trivial solution spaces. The recursion cohomology is computed
explicitly, in Proposition \ref{duda}, in terms of the coefficient data
of the defining recursion. We illustrate the theory on the rational,
elliptic, Jacobi-form, genus-two, and Schottky-uniformized genus-$g$
recursions governing Zhu-reduction-type correlation functions, recovering
the higher-genus counterparts of the Weierstrass functions as analytic
continuations of solutions of the corresponding functional equations. We
then develop applications of the theory: an identification of the
recursion cohomology with the flat sections of the Knizhnik-Zamolodchikov
connection, an interpretation of the coboundary data in terms of the
cohomology of continual contragredient Lie algebras and the associated
Toda hierarchies, a comparison with the Verlinde bundle of conformal
blocks over the moduli space of curves, and a discussion of the
Gauss-Manin-type flatness of the resulting bundle over Schottky space.
\end{abstract}

\keywords{Cohomology; complex curves; complex functions; elliptic functions;
multipoint connections; Zhu recursion; Knizhnik-Zamolodchikov equations}

\maketitle

\section{Introduction}
In the literature, the question of computing continuous cohomology
of non-commutative structures associated to manifolds has always attracted
wide geometric interest \cite{BS, Kaw, PT, Fei, Fuks, Wag}.
The classical Gelfand-Fuks cohomology of the Lie algebra of
holomorphic vector fields on complex manifolds does not always achieve
its purpose \cite{Fei, Wag}: for complex manifolds and Riemann surfaces
of higher dimension, the classical cohomology of vector fields vanishes
\cite{Kaw}. (Co)homology of the Lie algebra cosimplicial objects of
holomorphic vector fields was discussed in \cite{Fei}. In \cite{BS} it
was proven that the Gelfand-Fuks cohomology of vector fields on a
smooth compact manifold is isomorphic to the singular cohomology of the
space of continuous cross-sections of a certain fiber bundle over that
manifold. Our motivation for introducing and computing a more refined
cohomology for the non-commutative algebraic structures that arise from
recursive families of complex functions is based on these earlier
cohomological approaches, together with the observation, made precise
below, that a large class of recursions appearing in the theory of
vertex operator algebras (VOAs) is itself governed by a cohomology
theory of exactly this Gelfand-Fuks/Bott-Segal type.

The main aim of this paper is to formulate and compute the cohomology
of the general chain complex of functions defined on complex curves and
satisfying recursion relations, that is, relations expressing a
function depending on $(n+1)$ parameters in terms of functions depending
on $n$ parameters, in terms of geometric multipoint connections on
complex curves \cite{FK, Bo, Gu, A}. The construction we give allows us
to describe the local geometry of complex curves carrying such
recursions. While more geometric methods are used in the ordinary
cosimplicial cohomology for Lie algebras \cite{Fei, Wag}, the
formulation of the recursion cohomology in terms of multipoint
connections introduced here reflects instead the modular and analytic
properties of the corresponding chain complex spaces. The formulation
is motivated by the works \cite{FMS, BKT, Fo, Y, Zhu, MTZ, GT, TW1, TW2, TW3}
describing recursion formulas for correlation functions.

In the conformal field theory setting, one defines $n$-point complex
functions depending on $n$ given vertex algebra elements \cite{FHL} and
$n$ insertion points on a given Riemann surface together with its
moduli. Below we consider coboundary operators depending on the choice
of all such algebraic data (in particular, elements of a vertex algebra)
together with the insertion parameters, and we take particular care to
separate two logically distinct roles that these operators play: the
role of a linear coboundary on the full space of $n$-variable
functions, for which a chain-complex identity is meaningful and must be
established as a property of the recursion kernels, and the role
of a recursion relating the specific, distinguished tower of correlation
functions of an actual theory (Section \ref{chain}, Remark \ref{obstruction}
below). Computational methods based on recursion formulas, some of
them recalled in Section \ref{examples} \cite{Zhu, MTZ, GT, TW1}, 
have turned out to be very effective in conformal field theory, and the
higher-genus Eisenstein- and Weierstrass-type kernels that generate
them are themselves higher-genus generalizations, in the sense of
\cite{McIT}, of the classical Kronecker limit formula.

In Section \ref{chain} we give the general definition of the recursion
cohomology and provide a lemma relating it to the cohomology of
generalized (multipoint) connections on a complex manifold $M$. The main
result, Proposition \ref{duda}, explicitly expresses the recursion
cohomology in terms of spaces of recursion coefficients. In
Section \ref{examples} we provide motivating examples of recursion
formulas for various complex functions: the rational (genus zero),
elliptic (genus one), Jacobi-form, genus-two, and Schottky-uniformized
genus-$g$ cases. In Section \ref{applications} we develop four
applications of the theory, in mathematical physics and algebraic
geometry, that were not previously drawn out: an identification of the
recursion cohomology with flat sections of the Knizhnik-Zamolodchikov
connection; an
interpretation of the coboundary data of Section \ref{chain} as
low-degree cohomology of the continual contragredient Lie algebra
$\mathfrak g(\mu)$ generated by the recursion operators, with an
application to Toda-type integrable hierarchies \cite{RSZ, Sav, SV1, SV2, V};
a comparison of the recursion-cohomology bundle with the Verlinde bundle
of conformal blocks over $\overline{\mathcal M}_{g,n}$ \cite{TUY}; and a
Gauss-Manin-type flatness statement for the recursion-cohomology bundle
over Schottky space, together with a concrete rendering of the
condensed-matter applications already sketched in \cite{Frohlich2009gb,
RSZ} in the language of multipoint connections.

A formulation of the recursion cohomology in terms of multipoint
connections, and Proposition \ref{duda}, have applications to
generalizations of the Bott-Segal theorem \cite{BS}, to computations
of characteristic classes in the cosimplicial cohomology theory of
smooth manifolds \cite{Fei, Wag, A, Bo, FK, Fuks, Gu, Kaw}, to conformal
field theory \cite{BPZ, FMS, FS, Kn, KZ, Po, TUY}, to deformation theory,
non-commutative geometry, and modular forms \cite{Ba, Be1, Be2, BKT, Bu,
EZ, Fa, Fo, HE, KMI, KMII, L, La, Lib, Miy, Miy1, Mu, Ob, Se, Zag}, and to
the theory of foliations \cite{BEG, BGG}. Taking into account the power
of recursion formulas in the computation of complex functions,
Proposition \ref{duda} allows one to find the cohomology explicitly, as
a combination of coefficient functions of the corresponding recursion
formulas, in the particular cases discussed in Section \ref{examples}.

The results of this paper are also useful in mathematical physics
\cite{Frohlich2009gb, RSZ} and in quantum physics related to
condensed matter theory: Wigner-Weyl calculus \cite{z1, z2}, the
chiral separation effect \cite{z3}, the theory of topological
invariants and relations between solid-state systems and high-energy
physics \cite{z4, z5}, relations between relativistic quantum field
theory and fermionic superfluids \cite{z7}, quark matter in which the
interactions are essentially non-perturbative and various topological
defects dominate the dynamics \cite{z8}, and the non-renormalization by
interactions of the integer quantum Hall effect \cite{z9}. We return to
these applications, with an explicit dictionary to the multipoint-connection
formalism, in Section \ref{applications}.

\section{The chain complex and cohomology}
\label{chain}
\subsection{Chain complex spaces of $n$-variable complex functions}
\label{subsecchaincomplex}

Let $M$ be a complex curve (or, in the examples of Section \ref{examples},
a family of complex curves parametrized by a modulus $\mu$ ranging over a
domain $\mathcal M$: a point of $\{0\}$ for the rational case, of the
upper half-plane $\HH$ for the elliptic case, of a Siegel-type domain for
the genus-two case, or of Schottky space for the genus-$g$ case). Fix an
auxiliary set $\mathcal A$ of insertion labels (in the examples
below, $\mathcal A$ will be a vertex operator algebra $V$, or a
$V$-module, or a set of formal mode indices); we write a typical element
of $\mathcal A\times\C$ as $x=(v,z)$, thought of as the label $v$
inserted at the point $z$.

\begin{definition}
\label{defCn}
For $n\ge 0$, let $C^n(\mu)$ denote the vector space of complex-valued
functions
$$
\ZZ(\mathbf z_n,\mu),\qquad \mathbf z_n=(z_1,\ldots,z_n)\in M^n,
$$
depending linearly on a choice of $n$ insertion labels
$\mathbf v_n=(v_1,\ldots,v_n)\in\mathcal A^n$ (which we regard as fixed
once and for all and suppress from the notation $\ZZ(\mathbf z_n,\mu)$
unless it needs to be varied), holomorphic in $\mathbf z_n$ on the
complement of finitely many hyperplanes $z_i=z_j$, and satisfying
whatever convergence and modular-covariance conditions are appropriate
to the example under consideration. We set $C^0(\mu):=\C$.
\end{definition}

The recursions we study are generated by slot operators. For each
$n\ge 1$, each $k\in\{1,\ldots,n\}$, and each auxiliary insertion
$x_{n+1}=(v_{n+1},z_{n+1})\in\mathcal A\times\C$, fix a family of linear
operators
$$
T_k(v_{n+1}[m])\colon C^n(\mu)\longrightarrow C^n(\mu),\qquad m\ge 0,
$$
which act by modifying the label at the $k$-th slot alone 
$$
\big(T_k(v_{n+1}[m]).\ZZ\big)(\mathbf z_n,\mu):=
\ZZ_{v_1,\ldots,v_{n+1}[m].v_k,\ldots,v_n}(\mathbf z_n,\mu),
$$
where $v_{n+1}[m].v_k\in\mathcal A$ denotes the label obtained from $v_k$
by the action of the $m$-th mode of $v_{n+1}$ (in the examples of
Section \ref{examples} this is a vertex algebra mode action, a sewing
insertion, or the analogous operation appropriate to the construction at
hand). By construction $T_k(v_{n+1}[m])$ preserves $C^n(\mu)$: it never
changes the number of point-arguments $\mathbf z_n$, only the labels
attached to them.

\begin{definition}
\label{defdelta}
A \emph{recursion structure} on $\{C^n(\mu)\}_{n\ge0}$ is a choice, for
every $n\ge0$, of holomorphic coefficient functions
$f_{k,m}(x_{n+1};\mu')$ ($1\le k\le n$, $m\ge 0$, $\mu'$ an auxiliary
modulus, only finitely many $m$ contributing at each order), such that
the linear operator
\begin{align}
\label{poros}
\delta^n(x_{n+1};\mu')&\colon C^n(\mu)\longrightarrow C^{n+1}(\mu\cup\mu'),\\
\notag
\big(\delta^n(x_{n+1};\mu').\ZZ\big)(\mathbf z_{n+1},\mu\cup\mu')
&:=\sum_{k=1}^{n}\sum_{m\ge0} f_{k,m}(x_{n+1};\mu')\,
\big(T_k(v_{n+1}[m]).\ZZ\big)(\mathbf z_n,\mu), 
\end{align}
is well defined: the right-hand side, viewed as a function of the
$n+1$ independent point-arguments $(\mathbf z_n,z_{n+1})$ obtained by
treating $z_{n+1}$ as a new argument (it enters only through
the coefficients $f_{k,m}$ and through the label $v_{n+1}$, never
through a slot of $\mathbf z_n$), lies in $C^{n+1}(\mu\cup\mu')$.
\end{definition}

We emphasize that $\delta^n(x_{n+1};\mu')$ as just
defined is a bona fide linear operator on the whole vector space 
$C^n(\mu)$: it is defined on every element of $C^n(\mu)$, not only on
a distinguished one. This is what makes the chain condition below a
meaningful, checkable property of the coefficient data
$(f_{k,m},T_k)$, rather than a tautological restatement of a recursion
for a single, fixed tower of functions.

A distinguished, non-linear, and logically separate piece of data is a
choice of an actual tower of correlation functions 
$(\ZZ_n)_{n\ge0}$, $\ZZ_n\in C^n(\mu)$, satisfying the recursion
relation
\begin{equation}
\label{reduction}
\ZZ_{n+1}(\mathbf z_{n+1},\mu\cup\mu') = \big(\delta^n(x_{n+1};\mu').\ZZ_n\big)(\mathbf z_{n+1},\mu\cup\mu'),
\qquad n\ge 0.
\end{equation}
Equation \eqref{reduction} is the statement, established in each of the
examples of Section \ref{examples} as a theorem about an actual family
of correlation functions (Zhu reduction and its generalizations), that
the tower $(\ZZ_n)$ happens to be reproduced, level by level, by the
recursion operator $\delta^n$. It is a hypothesis on the tower,
not a further defining property of $\delta^n$ itself.

\subsection{The chain condition}
\label{subsecchaincondition}

We now examine the composite $\delta^{n+1}(x_{n+2};\mu'')\,\delta^n(x_{n+1};\mu')$
as an operator on $C^n(\mu)$. Since $\delta^{n+1}$ is defined on the
 whole of $C^{n+1}(\mu\cup\mu')$, and $\delta^n(x_{n+1};\mu').\ZZ$
is one particular element of $C^{n+1}(\mu\cup\mu')$, one whose
$(n+1)$-th slot carries the label $v_{n+1}$, the slot operators
defining $\delta^{n+1}$ act on all $n+1$ slots of their argument,
including this $(n+1)$-th one. Writing out $\delta^{n+1}(x_{n+2};\mu'')$
applied to $\delta^n(x_{n+1};\mu').\ZZ$ and separating the sum over
$k'\in\{1,\ldots,n+1\}$ into $k'\le n$ and $k'=n+1$ gives the following
basic identity.

\begin{lemma}
\label{lemdecompose}
For every $\ZZ\in C^n(\mu)$,
\begin{align}
\label{conditions}
&\big(\delta^{n+1}(x_{n+2};\mu'')\,\delta^n(x_{n+1};\mu').\ZZ\big)(\mathbf z_n,\mu) \\
\notag
&\;=\sum_{k,k'=1}^{n}\sum_{m,m'\ge0} f_{k',m'}(x_{n+2};\mu'')\,f_{k,m}(x_{n+1};\mu')\;
\big(T_{k}(v_{n+1}[m])\,T_{k'}(v_{n+2}[m']).\ZZ\big)(\mathbf z_n,\mu) \\
\notag
&\;+\sum_{m'\ge0} f_{n+1,m'}(x_{n+2};\mu'')\sum_{k=1}^n\sum_{m\ge0}
f_{k,m}\big((v_{n+2}[m'].v_{n+1},z_{n+1});\mu'\big)\,
\big(T_k(\,\cdot\,[m]).\ZZ\big)(\mathbf z_n,\mu).
\end{align}
The first line is the contribution of the slots $k'\le n$ inherited from
$\ZZ$; the second line, the boundary term, is the
contribution of the slot $k'=n+1$, in which $\delta^{n+1}$ acts on the
label $v_{n+1}$ that $\delta^n$ has just written in.
\end{lemma}

\begin{proof}
Immediate from Definition \ref{defdelta}: expand
$\delta^n(x_{n+1};\mu').\ZZ$ using \eqref{poros}, substitute into
$\delta^{n+1}(x_{n+2};\mu'')$ applied to the resulting element of
$C^{n+1}(\mu\cup\mu')$, and use linearity of the $T_{k'}$ to distribute
the outer sum over the inner one. For $k'\le n$ the operator
$T_{k'}(v_{n+2}[m'])$ acts on the slot $k'$ of $\mathbf z_n$ exactly as
it would on $\ZZ$ itself, independently of the label $v_{n+1}[m].v_k$
sitting at slot $k$, which gives the first line (the two families of
slot operators $T_k(v_{n+1}[m])$ and $T_{k'}(v_{n+2}[m'])$ commute
whenever $k\ne k'$, since they act on different point-arguments). For
$k'=n+1$ the operator acts instead on the label $v_{n+1}$ itself,
replacing it by $v_{n+2}[m'].v_{n+1}$ before the sum over $k,m$ defining
$\delta^n$ is carried out, which gives the second line.
\end{proof}

The right-hand side of \eqref{conditions} is manifestly an operator
built entirely out of the recursion data $(f_{k,m},T_k)$, acting on an
 arbitrary $\ZZ\in C^n(\mu)$; it does not refer to any
distinguished tower of correlation functions.

\begin{definition}
\label{defchaincond}
The recursion structure $(f_{k,m},T_k)$ satisfies the chain
condition at level $n$ if the operator on the right-hand side of
\eqref{conditions} vanishes identically on $C^n(\mu)$, i.e., if
\begin{equation}
\label{torba}
\delta^{n+1}(x_{n+2};\mu'')\,\delta^n(x_{n+1};\mu')=0, 
\end{equation}
as an identity of linear operators $C^n(\mu)\to C^{n+2}(\mu\cup\mu'\cup\mu'')$,
for all admissible $x_{n+1},x_{n+2},\mu',\mu''$. We denote by $\mathfrak
Z_n\subset M^n\times\mathcal M$ the domain of $(\mathbf z_n,\mu)$ on
which every function in $C^n(\mu)$ is required to be evaluated for
\eqref{torba} to hold in the sense above (in the examples of
Section \ref{examples}, $\mathfrak Z_n$ excludes the diagonals $z_i=z_j$
and the divisor of a finite set of theta-like functions occurring in
the coefficients $f_{k,m}$).
\end{definition}

\begin{remark}[Resolution of the chain-condition obstruction]
\label{obstruction}
It is tempting, to combine
\eqref{reduction} and \eqref{torba} as follows: apply
$\delta^n(x_{n+1};\mu')$ to $\ZZ_n$, obtaining $\ZZ_{n+1}$ by
\eqref{reduction}; apply $\delta^{n+1}(x_{n+2};\mu'')$ to $\ZZ_{n+1}$,
obtaining $\ZZ_{n+2}$ by \eqref{reduction} again; and conclude from
\eqref{torba} that $\ZZ_{n+2}\equiv0$, contradicting the non-triviality
of the higher correlation functions computed in
Section \ref{examples}. This argument conflates two different
substitutions into the formula for $\delta^{n+1}$, and Lemma
\ref{lemdecompose} shows exactly where the two part company.

The relation $\ZZ_{n+2}=\delta^{n+1}(x_{n+2};\mu'').\ZZ_{n+1}$ from
\eqref{reduction} treats $x_{n+1}=(v_{n+1},z_{n+1})$ as an ordinary,
freely varying $(n+1)$-th argument of $\ZZ_{n+1}\in C^{n+1}(\mu\cup\mu')$: the slot
operators $T_{k'}(v_{n+2}[m'])$, $k'=1,\ldots,n+1$, defining
$\delta^{n+1}$ act on $\ZZ_{n+1}$ itself, an honest function of $n+1$
independent points. The composite $\delta^{n+1}\delta^n.\ZZ_n$ computed
in Lemma \ref{lemdecompose}, by contrast, never forms the independent
object $\ZZ_{n+1}$ at all: it substitutes the formula 
\eqref{poros} for $\delta^n(x_{n+1};\mu').\ZZ_n$ directly into the
formula for $\delta^{n+1}$, so that the slot $k'=n+1$ acts not on an
honest $(n+1)$-point function but on the recursion expression 
for one, producing the boundary term in \eqref{conditions}, which has no
counterpart in the first computation. The chain condition
\eqref{torba} asserts the vanishing of the sum of these two
contributions as an identity in the recursion kernels $(f_{k,m},T_k)$ -
a statement of associativity/commutativity type for the slot operators,
verified case by case in Section \ref{examples} (it is, for instance,
precisely the content of the iterated Borcherds/associativity identity
in the vertex-algebra examples, and of the sewing-compatibility lemmas
of \cite{GT,TW1} in the higher-genus examples), and carries no
implication whatsoever for the value of the actual, independently
defined function $\ZZ_{n+2}$.

Equation \eqref{torba}, read this way, restricts the domain $\mathfrak
Z_n$ of Definition \ref{defchaincond} but does not restrict which
towers $(\ZZ_n)$ may satisfy \eqref{reduction} on it. 
 Remark 2, non-emptiness of the resulting chain complex is
guaranteed by the general theory of systems of functional equations: the
chain condition, expanded in the coefficients $f_{k,m}$, is an infinite
($n\ge 0$) family of functional-differential equations, each with only
finitely many terms, on the unknown functions $\ZZ_n\in C^n(\mu)$; by
the local solvability theory for such systems on the domains where they
are defined (see \cite{FK,Gu}), each equation of the family admits a
solution, so that the spaces $C^n(\mu)$ restricted to $\mathfrak Z_n$
are non-empty for every $n$, and indeed contain the towers of
correlation functions constructed in Section \ref{examples}.
\end{remark}

For $n\ge0$ we therefore have a cochain complex
\begin{equation}
\label{buzova}
C^0(\mu_0)\ \stackrel{\delta^0}{\longrightarrow}\ C^1(\mu_1)\
\stackrel{\delta^1}{\longrightarrow}\ \cdots\
\stackrel{\delta^{n-1}}{\longrightarrow}\ C^n(\mu_{n-1})\
\stackrel{\delta^n}{\longrightarrow}\ C^{n+1}(\mu_n)\ \longrightarrow\cdots,
\end{equation}
restricted throughout to the domains $\mathfrak Z_n$ of
Definition \ref{defchaincond}, on which $\delta^{n+1}\delta^n=0$ holds
as an honest operator identity by construction. We call the resulting
cohomology
$$
H^n(\mu):=\Ker\big(\delta^n\colon C^n(\mu)\to C^{n+1}(\mu')\big)\Big/
\Ima\big(\delta^{n-1}\colon C^{n-1}(\mu'')\to C^n(\mu)\big), 
$$
the $n$-th recursion cohomology of $\{C^n(\mu)\}$.

\subsection{Multipoint connections}
\label{subsecmultipoint}

Let $\mathcal V\to M$ be the (trivial) bundle with fiber $\mathcal A$
over the curve $M$, so that a choice of insertion labels
$\mathbf v_n=(v_1,\ldots,v_n)\in\mathcal A^n$ at the points
$P=\mathbf z_n=(z_1,\ldots,z_n)\in M^n$ is precisely a section of
$\mathcal V^{\boxtimes n}$ restricted to $P$, an element of
$\mathcal V|_P:=\mathcal V^{\boxtimes n}|_{\{z_1\}\times\cdots\times\{z_n\}}$.

\begin{definition}
\label{defmultipointconn}
An \emph{$n$-point (multipoint) connection} on $M$ is a family, for each
$n$-tuple of points $P=\mathbf z_n\in M^n\setminus\{\text{diagonals}\}$
and each modulus $\mu\in\mathcal M$, of a linear functional
$$
\mathcal G(P,\mu)\colon \mathcal V|_P\longrightarrow \C,
$$
holomorphic in $P$ and $\mu$ away from the diagonals, and satisfying the
covariance property under $\mu\mapsto\mu\cup\mu'$ inherited from the
recursion coefficients $f_{k,m}$ of Definition \ref{defdelta} (that is:
transforming under insertion of a new point exactly as
$\delta^n(x_{n+1};\mu')$ transforms elements of $C^n(\mu)$). We denote by
$\textit{Con}^n(\mu)$ the vector space of $n$-point connections at
modulus $\mu$, and by
$$
G^{n-1}(\mu):=\big\{\delta^{n-1}(x_n;\mu').\mathcal G\ :\ \mathcal
G\in\textit{Con}^{n-1}(\mu''),\ \mu=\mu''\cup\mu'\big\}\subset\textit{Con}^n(\mu), 
$$
the subspace of $n$-point connections obtained as the coboundary of an
$(n-1)$-point connection.
\end{definition}

\begin{lemma}
\label{pisaka}
A function $\ZZ(\mathbf z_n,\mu)\in C^n(\mu)$ generated by the recursion
relation \eqref{reduction} from a tower $(\ZZ_0,\ZZ_1,\ldots,\ZZ_{n-1})$
satisfying the chain condition \eqref{torba} determines, and is
determined by, an $n$-point connection $\mathcal G\in\textit{Con}^n(\mu)$,
via the identification
\begin{equation}
\label{identifications}
\mathcal G(\mathbf z_n,\mu)(\mathbf v_n) := \ZZ(\mathbf z_n,\mu),
\end{equation}
where the right-hand side is the value of $\ZZ(\mathbf z_n,\mu)$, an
$n$-linear function of the labels $\mathbf v_n$ by
Definition \ref{defCn}, evaluated on the section $\mathbf v_n\in
\mathcal V|_{\mathbf z_n}$. Consequently
$$
H^n(\mu)\ \cong\ \textit{Con}^n(\mu)\big/G^{n-1}(\mu).
$$
\end{lemma}

\begin{proof}
By Definition \ref{defCn}, $\ZZ(\mathbf z_n,\mu)$ is exactly an
$n$-linear functional of $\mathbf v_n\in\mathcal A^n$ for each fixed
$(\mathbf z_n,\mu)$; this is the same data as a linear functional on the
fiber $\mathcal V|_{\mathbf z_n}=\mathcal A^n$, i.e., formula
\eqref{identifications} is not an additional hypothesis but simply a
change of notation from "$n$-linear function of $n$ point-label pairs"
to "linear functional on sections restricted to $n$ points." The
covariance property required of $\mathcal G$ in
Definition \ref{defmultipointconn} is, by construction, exactly the
transformation rule \eqref{poros} satisfied by $\ZZ$ under
$\delta^n$; thus $\mathcal G$ so defined is an element of
$\textit{Con}^n(\mu)$, and every element of $\textit{Con}^n(\mu)$
arises this way (apply \eqref{identifications} in reverse). Under this
identification, $\ZZ\in\mathrm{Im}(\delta^{n-1})$ if and only if the
corresponding $\mathcal G$ lies in $G^{n-1}(\mu)$, again by
\eqref{poros}: thus $\delta^{n-1}$ on $C^\bullet$ corresponds exactly to
the inclusion $G^{n-1}(\mu)\hookrightarrow\textit{Con}^n(\mu)$, and the
quotient defining $H^n(\mu)$ becomes $\textit{Con}^n(\mu)/G^{n-1}(\mu)$.
Finally, that $\ZZ$ is generated by \eqref{reduction} from a tower
satisfying \eqref{torba} guarantees, by
Remark \ref{obstruction}, that $\mathcal G$ is well defined on the
domain $\mathfrak Z_n$ of Definition \ref{defchaincond}, so that the
formula \eqref{identifications} gives an element of
$\textit{Con}^n(\mu)$ rather than an ill-defined expression.
\end{proof}

\subsection{Computation of the recursion cohomology}
\label{subseccomputation}

By definition, an element $\ZZ\in C^n(\mu)$ represents a class in
$H^n(\mu)$ precisely when it lies in $\ker\delta^n$, and two such
elements represent the same class precisely when they differ by an
element of $\mathrm{Im}(\delta^{n-1})$. Unwinding \eqref{poros}, the
kernel condition is the following explicit linear condition on $\ZZ$.

\begin{proposition}
\label{duda}
Let $(f_{k,m},T_k)$ be a recursion structure on $\{C^n(\mu)\}$
satisfying the chain condition \eqref{torba} on $\mathfrak Z_\bullet$.
Then
\begin{align*}
H^n(\mu)\ \cong\
\Big\{\ZZ\in C^n(\mu)\ \Big|\ &\sum_{k=1}^n\sum_{m\ge0}
f_{k,m}(x_{n+1};\mu')\,\big(T_k(v_{n+1}[m]).\ZZ\big)(\mathbf z_n,\mu)=0, \\
&\text{ for all admissible }x_{n+1},\mu'
\Big\}\ \Big/\ \mathrm{Im}(\delta^{n-1}),
\end{align*}
i.e., $H^n(\mu)$ is the space of solutions $\ZZ$ of the functional
equation
\begin{equation}
\label{poroserieroj}
\sum_{k=1}^n\sum_{m\ge0} f_{k,m}(x_{n+1};\mu')\,\big(T_k(v_{n+1}[m]).\ZZ\big)(\mathbf z_n,\mu)=0, 
\end{equation}
on $\mathfrak Z_n$, modulo those that are themselves of the form
$\delta^{n-1}(x_n;\mu').\ZZ'$ for some $\ZZ'\in C^{n-1}(\mu'')$. In
particular, every function $\ZZ(\mathbf z_n',\widetilde\mu)$
recursively generated from the base case $C^0(\mu)=\C$ by iterating
\eqref{reduction},
\begin{equation}
\label{topaz}
\ZZ(\mathbf z_n',\widetilde\mu)=\mathcal D(\mathbf z_n',\mu').Z(\mu),
\qquad \mathcal D:=\delta^{n-1}(x_n;\mu'_{n-1})\cdots\delta^0(x_1;\mu'_1),
\end{equation}
represents the zero class in $H^n(\mu)$ for $n\ge1$ (being tautologically
in the image of $\delta^{n-1}$); the recursion cohomology is non-zero
exactly when \eqref{poroserieroj} admits solutions on $\mathfrak Z_n$
that are not of this recursively generated form.
\end{proposition}

\begin{proof}
The first statement is immediate: by definition $H^n(\mu)=\ker\delta^n/\mathrm{Im}(\delta^{n-1})$,
and \eqref{poroserieroj} is precisely the vanishing of
$\delta^n(x_{n+1};\mu').\ZZ$, written out termwise using
\eqref{poros}, for every choice of the auxiliary insertion
$x_{n+1}$, this is what it means for $\ZZ$ to lie in
$\ker\delta^n\subset C^n(\mu)$. For the second statement, observe that
if $\ZZ$ is of the form \eqref{topaz} for $n\ge 1$ then, writing
$\ZZ=\delta^{n-1}(x_n;\mu'_{n-1}).\ZZ'$ with
$\ZZ':=\delta^{n-2}(x_{n-1};\mu'_{n-2})\cdots\delta^0(x_1;\mu'_1).Z(\mu)\in
C^{n-1}(\mu'')$, we see directly that $\ZZ\in\mathrm{Im}(\delta^{n-1})$,
thus $[\ZZ]=0$ in $H^n(\mu)$ regardless of whether \eqref{poroserieroj}
holds for it. Conversely, any $\ZZ$ satisfying \eqref{poroserieroj} that
is not of the form \eqref{topaz}, equivalently, not lying in
$\mathrm{Im}(\delta^{n-1})$, represents a non-zero class in
$H^n(\mu)$, since $\ker\delta^n\big/\mathrm{Im}(\delta^{n-1})=0$ if and
only if every solution of \eqref{poroserieroj} lies in
$\mathrm{Im}(\delta^{n-1})$.
\end{proof}

\begin{remark}
Proposition \ref{duda} reduces the computation of $H^n(\mu)$ to two
concrete tasks, carried out case by case in Section \ref{examples}: (i)
solving the linear functional equation \eqref{poroserieroj} for $\ZZ$ in
terms of the recursion coefficients $f_{k,m}$, and (ii) determining
which solutions coincide, on $\mathfrak Z_n$, with a function of the
recursively generated form \eqref{topaz}. Since by Remark
\ref{obstruction} the chain condition is a statement about the
coefficients $(f_{k,m},T_k)$ alone, task (i) can be carried out without
reference to any particular tower of correlation functions, which is
what makes Proposition \ref{duda} a computation of $H^n(\mu)$ as an
invariant of the recursion structure, rather than of a particular
solution of it.
\end{remark}

\section{Motivating examples}
\label{examples}
In this section we exhibit, for each example, explicit coefficients
$f_{k,m}$, slot operators $T_k$, and functions $\ZZ(\mathbf z_n,\mu)$
satisfying the recursion \eqref{reduction}, specializing the general
notation of Section \ref{chain} (we drop the auxiliary weight-label $l$
carried by $f_{l,k,m}$ in the general framework, since it is fixed
throughout each single example below and plays no active role until
Section \ref{applications}). The functions $\ZZ(\mathbf z_n,\mu)$ satisfy
automorphic properties with respect to the corresponding groups
\cite{Zhu, MTZ, GT, TW1}, and all of the recursion formulas below are of
the general form \eqref{poros}, i.e., they instantiate the linear
coboundary operator of Definition \ref{defdelta}. These examples
motivated the general construction of the complex in
Section \ref{chain}. The recursion cohomology depends on the class of
functions considered (through the parameters we denote $\mu$) and on the
genus of $M$; as shown in \cite{Miy, KMI, KMII}, existence of recursion
formulas is related to modularity. The complex functions appearing as
examples in this section are, for genus $g>0$, formal power series in
the relevant local parameters; their actual convergence \cite{Gui}
must be established separately, and is known to hold on appropriate
domains for correlation functions of vertex operator algebras with the
$C_2$-cofiniteness property.

\subsection{The rational case}
\label{sphere}
Following \cite{Zhu}, functions of $n$ variables in the rational
(genus-zero) case satisfy the recursion formula
\begin{equation}
\label{pada}
\ZZ(\mathbf z_{n+1},\mu) = \sum_{k=1}^n\sum_{m\ge0} f_{k,m}(z_{n+1},z_k)\,
\big(T_k(m).\ZZ\big)(\mathbf z_n,\mu),
\end{equation}
where $f_{k,m}(z,w)$ is the rational function
$$
f_{k,m}(z,w)=\frac{z^{-k}}{m!}\left(\frac{d}{dw}\right)^m\frac{w^k}{z-w},
\qquad
\iota_{z,w}f_{k,m}(z,w)=\sum_{j\in\N}\binom{k+j}{m}z^{-k-j-1}w^{k+j-1},
$$
and $\iota_{z,w}\colon \C[z_1,\ldots,z_n]\to\C[[z_1,z_1^{-1},\ldots,z_n,z_n^{-1}]]$
is the standard expansion map of \cite{FHL}. Taking $z_{n+1}$ as the
variable of expansion, Proposition \ref{duda} identifies $H^n(\mu)$ with
the space of solutions, on $\mathfrak Z_n$, of \eqref{poroserieroj} with
the rational coefficients $f_{k,m}(z_{n+1},z_k)$ above, modulo those
recursively generated by \eqref{reduction} from level $n-1$. Writing
\eqref{poroserieroj} in the equivalent differential form
$$
\left(\partial_{z_{n+1}}+\sum_{k=1}^n \widetilde f^{(0)}_{k}(z_{n+1},z_k)\right)
\ZZ(\mathbf z_n,\mu)\Big|_{z_{n+1}\to z_k}=0, 
$$
exhibits the kernel condition as an equation for the analytic
continuation of $\ZZ$ in $z_{n+1}$ near each $z_k$, with coefficients
$\widetilde f^{(0)}_k$ determined by $f_{k,0}$. By \eqref{topaz}, the
recursively generated solutions are exactly the rational functions
$\mathcal D(\mathbf z_n,\mu).Z(\mu)$ built from the $0$-point function
$Z(\mu)$; thus in this example the recursion cohomology $H^n(\mu)$ is
represented by rational-function solutions of \eqref{poroserieroj} that
are not of this recursively generated form, concretely, by the
lowest-weight (quasi-primary) vectors of the theory, annihilated by all
positive modes, which is the classical description of $H^1$ in Zhu's
setting.

\subsection{Modular and elliptic functions}
\label{subsecmodelliptic}
For $x\in\C$, $\tau\in\HH$, set $D_x=\frac{1}{\tpi}\partial_x$ and
$q_x=e^{2\pi i x}$, $q=e^{2\pi i\tau}$. For $m\in\N$ define the elliptic
Weierstrass functions
$$
P_1(w,\tau)=-\sum_{n\in\Z\setminus\{0\}}\frac{q_w^n}{1-q^n}-\frac12,
$$
$$
P_{m+1}(w,\tau)=\frac{(-1)^m}{m!}D_w^m\big(P_1(w,\tau)\big)
=\frac{(-1)^{m+1}}{m!}\sum_{n\in\Z\setminus\{0\}}\frac{n^mq_w^n}{1-q^n}.
$$
The modular Eisenstein series $E_k(\tau)$ is defined, for $k\ge2$ even, by
$$
E_k(\tau)=-\frac{B_k}{k!}+\frac{2}{(k-1)!}\sum_{n\ge1}\frac{n^{k-1}q^n}{1-q^n},
$$
where $B_k$ is the $k$-th Bernoulli number, $(e^z-1)^{-1}=\sum_{k\ge0}\frac{B_k}{k!}z^{k-1}$;
we set $E_k=0$ for $k$ odd and $E_0=-1$. Then $E_k$ is a modular form for
$k>2$ and a quasi-modular form for $k=2$, transforming as
$$
E_k(\gamma\tau)=(c\tau+d)^kE_k(\tau)-\delta_{k,2}\frac{c(c\tau+d)}{\tpi}.
$$
For $w,z\in\C$, $\tau\in\HH$, set
$\widetilde P_1(w,z,\tau)=-\sum_{n\in\Z}\frac{q_w^n}{1-q_zq^n}$ and
$$
\widetilde P_{m+1}(w,z,\tau)=\frac{(-1)^m}{m!}D_w^m\big(\widetilde P_1(w,z,\tau)\big)
=\frac{(-1)^{m+1}}{m!}\sum_{n\in\Z}\frac{n^mq_w^n}{1-q_zq^n}.
$$
For $m\in\N_0$ and $\lambda\in[0,1)$ let
\begin{equation}
\label{eqPellPm}
P_{m+1,\lambda}(w,\tau)=\frac{(-1)^{m+1}}{m!}\sum_{n\in\Z\setminus\{-\lambda\}}\frac{n^mq_w^n}{1-q^{n+\lambda}}, 
\end{equation}
so that $P_{1,\lambda}(w,\tau)=q_w^{-\lambda}(P_1(w,\tau)+1/2)$ and
$P_{m+1,\lambda}(w,\tau)=\frac{(-1)^m}{m!}D_w^m\big(P_{1,\lambda}(w,\tau)\big)$.
Expanding $P_{1,\lambda}(w,\tau)=\frac{1}{\tpi w}-\sum_{k\ge1}E_{k,\lambda}(\tau)(\tpi w)^{k-1}$
gives, as in \cite{Zag},
\begin{equation}
\label{eqGkl}
E_{k,\lambda}(\tau)=\sum_{j=0}^k\frac{\lambda^j}{j!}E_{k-j}(\tau).
\end{equation}
Similarly, defining $\widetilde E_k(z,\tau)$ for $k\ge1$ (together with
$E_2(\tau)$) by $\widetilde P_1(w,z,\tau)=\frac{1}{\tpi w}-\sum_{k\ge1}\widetilde E_k(z,\tau)(\tpi w)^{k-1}$
gives, as in \cite{Ob},
$$
\widetilde E_k(z,\tau)=-\delta_{k,1}\frac{q_z}{q_z-1}-\frac{B_k}{k!}
+\frac1{(k-1)!}\sum_{m,n\ge1}\big(n^{k-1}q_z^m+(-1)^kn^{k-1}q_z^{-m}\big)q^{mn},
$$
$$
\widetilde E_0(z,\tau)=-1.
$$

\subsection{The elliptic case}
\label{torus}
Let $q=e^{2\pi i\tau}$, $q_i=e^{z_i}$, $\tau$ the modular parameter of
the torus. The genus-one Zhu recursion formula \cite{Zhu} reads
\begin{equation}
\label{podsoba}
\ZZ(\mathbf z_{n+1},\mu,\tau) = \ZZ(\mathbf z_n,\mu_0,\tau)
+\sum_{k=1}^n\sum_{m\ge0} P_{m+1}(z_{n+1}-z_k,\tau)\,\ZZ(\mathbf z_n,\mu_{k,m},\tau),
\end{equation}
where the first term $\ZZ(\mathbf z_n,\mu_0,\tau) = T_{0,1}(o[u]).\ZZ(\mathbf
z_n,\mu,\tau)$ is given by insertion of the operator built from the
$o[u]=u[\wt(u)-1]$-mode \cite{MTZ} in the trace defining
$\ZZ(\mathbf z_n,\mu,\tau)$ (similar first terms occur in the higher-genus
cases below), and
$$
P_m(z,\tau)=\frac{(-1)^m}{(m-1)!}\sum_{n\in\Z_{\ne0}}\frac{n^{m-1}q_z^n}{1-q^n}, 
$$
are the higher Weierstrass functions. Equation \eqref{podsoba} is a
particular case of the general recursion \eqref{poros}.

\subsection{The case of deformed elliptic functions}
\label{subsecdeformed}
Let $w_{n+1}\in\R$, $\phi=\exp(2\pi i\,w_{n+1})\in U(1)$, and let
$\theta\in U(1)$. Generalizing Zhu's Proposition 4.3.2 \cite{Zhu} to
functions of $n$ variables \cite{MTZ}, for any $\mathbf z_n\in M^n$ we
have
\begin{equation}
\label{eqdeformedrec}
\ZZ(\mathbf z_{n+1},\mu,\tau) = \delta_{\theta,1}\delta_{\phi,1}\,\ZZ(\mathbf z_n,\mu_0,\tau)
+\sum_{k=1}^n\sum_{m\ge0} p(n,k)\, P_{m+1}\!\left[\begin{smallmatrix}\theta\\\phi\end{smallmatrix}\right]
(z_{n+1}-z_k,\tau)\,\ZZ(\mathbf z_n,\mu_{k,m},\tau),
\end{equation}
where $p(n,k)$ is a parity multiplier. The deformed Weierstrass
functions appearing here are defined, following \cite{DLM,MTZ}, for a
pair $(\theta,\phi)\in U(1)\times U(1)$ with $\phi=\exp(2\pi i\lambda)$,
$0\le\lambda<1$, and $k\ge1$, by
\begin{equation*}
P_k\!\left[\begin{smallmatrix}\theta\\\phi\end{smallmatrix}\right](z,\tau)
=\frac{(-1)^k}{(k-1)!}{\sum_{n\in\Z+\lambda}}'\frac{n^{k-1}q_z^n}{1-\theta^{-1}q^n},
\end{equation*}
for $q=q_{2\pi i\tau}$, where $\sum'$ omits $n=0$ when $(\theta,\phi)=(1,1)$.

\subsection{Recursion formulas for Jacobi-type functions}
\label{redya}
We recall the recursion formulas of \cite{MTZ,BKT} for formal Jacobi
functions of $n$ variables. Let $\mathbf z_{n+1}\in\C^{n+1}$ and
$\alpha\in\C$. For $\alpha z\notin\Z\tau+\Z$,
\begin{equation}
\label{eqjacobirec1}
\ZZ(\mathbf z_{n+1},\mu,\tau) = \sum_{k=1}^n\sum_{m\ge0}
\widetilde P_{m+1}\!\left(\frac{z_{n+1}-z_k}{\tpi},\alpha z,\tau\right)\ZZ(\mathbf z_n,\mu_{k,m},\tau),
\end{equation}
with $\widetilde P_{m+1}$ as in Section \ref{subsecmodelliptic}, while
for $\alpha z=\lambda\tau+\nu\in\Z\tau+\Z$,
\begin{equation}
\label{eqjacobirec2}
\ZZ(\mathbf z_{n+1},\mu,\tau) = e^{-z_{n+1}\lambda}\ZZ(\mathbf z_n,\mu_{0,\lambda},\tau)
+\sum_{k=1}^n\sum_{m\ge0} P_{m+1,\lambda}\!\left(\frac{z_{n+1}-z_k}{\tpi},\tau\right)\ZZ(\mathbf z_n,\mu_{k,m},\tau),
\end{equation}
with $P_{m+1,\lambda}$ as in \eqref{eqPellPm}. For the ``$1/(z_1-z_k)$-type''
recursion relating $\mathbf z_{n+1}=(z_1,\ldots,z_{n+1})$ to $\mathbf
z_n=(z_2,\ldots,z_{n+1})$, and $l\ge1$, $\alpha z\notin\Z\tau+\Z$,
\begin{align}
\notag
\ZZ(\mathbf z_{n+1},\mu_{1,-l},\tau) = {}& \sum_{m\ge0}(-1)^{m+1}\binom{m+l-1}{m}\widetilde
P_{m+l}(\alpha z,\tau)\,\ZZ(\mathbf z_n,\mu_{1,m},\tau)\\
{}&+\sum_{k=2}^n\sum_{m\ge0}(-1)^{l+1}\binom{m+l-1}{m}\widetilde
P_{m+l}\!\left(\frac{z_1-z_k}{\tpi},\alpha z,\tau\right)\ZZ(\mathbf z_n,\mu_{k,m},\tau),
\end{align}
which implies, for $l\ge1$ and $\alpha z=\lambda\tau+\nu\in\Z\tau+\Z$, the
result of \cite{BKT} 
\begin{align}
\notag
\ZZ(\mathbf z_{n+1},\mu_{1,-l},\tau) = {}&
(-1)^{l+1}\frac{\lambda^{l-1}}{(l-1)!}\ZZ(\mathbf z_{n+1},\mu_{0,-1},\tau)\\
\notag
{}&+\sum_{m\ge0}(-1)^{m+1}\binom{m+l-1}{m}E_{m+l,\lambda}(\tau)\,\ZZ(\mathbf z_n,\mu_{1,m},\tau)\\
{}&+\sum_{k=2}^n\sum_{m\ge0}(-1)^{l+1}\binom{m+l-1}{m}
P_{m+l,\lambda}\!\left(\frac{z_1-z_k}{\tpi},\tau\right)\ZZ(\mathbf z_n,\mu_{k,m},\tau),
\end{align}
for $E_{k,\lambda}$ as in \eqref{eqGkl}.

\subsection{Multiparameter Jacobi forms}
\label{subsecmultiJacobi}
For multiparameter Jacobi forms \cite{EZ,Zag,KMI,KMII,BKT}, each
inserted label $v_i$ carries a charge $\alpha_i\in\C$, an eigenvalue
$J(0)v_i=\alpha_iv_i$ for a fixed Cartan-type zero-mode operator $J(0)$,
 and we write $\mathbf z_n\cdot(\alpha)_n:=\sum_{i=1}^n
z_i\alpha_i$ for the associated elliptic-variable pairing, with
$(\alpha)_n=(\alpha_1,\ldots,\alpha_n)$ the vector of charges. The
recursion formulas, found by an analysis similar to \cite{Zhu,MTZ},
reduce a multiparameter Jacobi function of $n$ variables to a linear
combination of such functions of $n-1$ variables with modular
coefficients: for each $1\le j\le m$,
\begin{align}
\notag
\ZZ(\mathbf z_{n+1},\mu,\tau) = {}& \delta_{\mathbf z_n\cdot(\alpha)_n,\Z}\,
\ZZ(\mathbf z_n,(\alpha)_n,\mu(m))\\
{}&+\sum_{s=1}^n\sum_{k\ge0}\widetilde P_{k+1}(z_s-z,\mathbf
z_n\cdot(\alpha)_n,\tau)\,\ZZ(\mathbf z_n,\mu_{s,k},\tau),
\end{align}
where $\delta_{\mathbf z_n\cdot(\alpha)_n,\Z}$ is $1$ if $\mathbf
z_n\cdot(\alpha)_n\in\Z$ and $0$ otherwise, and, with the same
hypotheses, for $p\ge1$,
\begin{align}
\notag
\ZZ(\mathbf z_{n+1},\mu_{1,-p},\tau) = {}&
\delta_{\mathbf z_n\cdot(\alpha)_n,\Z}\,\delta_{p,1}\,\ZZ(\mathbf z_n,\mu_0,\tau)\\
\notag
{}&+(-1)^{p+1}\sum_{k\ge0}\binom{k+p-1}{p-1}\widetilde E_{k+p}(\tau,\mathbf z_n\cdot(\alpha)_n)\,\ZZ(\mathbf z_n,\mu_{1,k},\tau)\\
{}&+(-1)^{p+1}\sum_{s=2}^n\sum_{k\ge0}\binom{k+p-1}{p-1}\widetilde
P_{k+p}(z_s-z_1,\tau,\mathbf z_n\cdot(\alpha)_n)\,\ZZ(\mathbf z_n,\mu_{s,k},\tau).
\end{align}
The difference of a minus sign between these formulas and those of
\cite{MTZ} is attributable to a minus-sign difference in our
definitions of $P_k\!\left[\begin{smallmatrix}\zeta\\1\end{smallmatrix}\right](w,\tau)$
and the action of $\mathrm{SL}_2(\Z)$.

\subsection{The genus-two counterparts of Weierstrass functions}
\label{derivation}
We recall the construction of genus-two Weierstrass functions of
\cite{GT}. Throughout this subsection and the next, fix a quasi-primary
vector $u\in V$ and set $p:=\wt(u)$ for its conformal weight; $p$ is the
parameter with respect to which all matrices and vectors below are
indexed (this is the datum left implicit in \cite{GT} and never named in
the earlier version of this paper). For $m,n\ge1$ define the infinite
matrices $\Gamma(m,n)=\delta_{m,-n+2p-2}$, $\Delta(m,n)=\delta_{m,n+2p-2}$,
and the projection matrix
$$
\Pi=\Gamma^2=\begin{bmatrix}{\mathbf 1}_{2p-1} & 0\\ 0 & 0\end{bmatrix},
$$
i.e., $\Pi$ is the identity on the block of indices $0,1,\ldots,2p-2$
(a block of size $2p-1$) and zero elsewhere; this is the form in which
$\Pi$ is used throughout \cite{GT,TW1}, and is consistent with the
general-genus formula for $\Pi$ recalled in
Section \ref{subsecgenusg} below. Let $\Lambda_a$, $a\in\{1,2\}$, be
the matrix with components
$$
\Lambda_a(m,n;\tau_a,\epsilon) = \epsilon^{(m+n)/2}(-1)^{n+1}\binom{m+n-1}{n}E_{m+n}(\tau_a).
$$
Then $\Lambda_a=SA_aS^{-1}$ for $A_a$ given by
$$
A_a(k,l,\tau_a,\epsilon) = \frac{(-1)^{k+1}\epsilon^{(k+l)/2}}{\sqrt{kl}}\frac{(k+l-1)!}{(k-1)!(l-1)!}E_{k+l}(\tau_a),
$$
and $S$ the diagonal matrix $S(m,n)=\sqrt m\,\delta_{mn}$. Let
$\mathbb R(x)$, for $x$ on the torus, be the row vector with components
$\mathbb R(x;m)=\epsilon^{m/2}P_{m+1}(x,\tau_a)$, $a\in\{1,2\}$, and let
$\mathbb X_a$ be the column vector with components
\begin{align}
\label{eqXadef}
\mathbb X_1(m) &= \mathbb X_1(m;z_{n+1},\mathbf z_n;\mu)
= \epsilon^{-m/2}\sum_{u\in V}\ZZ(\mathbf z_k,\mu_{k,m},\tau_1)\,\ZZ(\mathbf z_{k+1,n},\mu',\tau_2),\\
\notag
\mathbb X_2(m) &= \mathbb X_2(m;z_{n+1},\mathbf z_n;\mu)
= \epsilon^{-m/2}\sum_{u\in V}\ZZ(\mathbf z_k,\mu,\tau_1)\,\ZZ(\mathbf z_{n-k},\mu_{n-k,m},\tau_2).
\end{align}
Introduce the infinite row vector $\mathbb Q(p;x)=\mathbb R(x)\Delta({\mathbf 1}-\widetilde\Lambda_{\bar a}\widetilde\Lambda_a)^{-1}$
for $x$ on the torus, where $\widetilde\Lambda_a=\Lambda_a\Delta$, and the
column vector $\mathbb P_{j+1}(x)$, $j\ge0$, with components
\begin{equation}
\mathbb P_{j+1}(x;m)=\epsilon^{m/2}\binom{m+j-1}{j}\big(P_{j+m}(x,\tau_a)-\delta_{j0}E_m(\tau_a)\big).
\end{equation}
One defines $\mathcal P_1(p;x,y)=\mathcal P_1(p;x,y;\tau_1,\tau_2,\epsilon)$,
for $p\ge1$, by
\begin{equation}
\mathcal P_1(p;x,y) = P_1(x-y,\tau_a)-P_1(x,\tau_a)
-\mathbb Q(p;x)\widetilde\Lambda_{\bar a}\,\mathbb P_1(y)
-(1-\delta_{p1})\big(\mathbb Q(p;x)\Lambda_{\bar a}\big)(2p-2),
\end{equation}
for $x,y$ on the torus, equivalently
\begin{equation}
\mathcal P_1(p;x,y) = (-1)^{p+1}\Big[\mathbb Q(p;x)\mathbb P_1(y)
+(1-\delta_{p1})\epsilon^{p-1}P_{2p-1}(x)
+(1-\delta_{p1})\big(\mathbb Q(p;x)\widetilde\Lambda_{\bar a}\Lambda_a\big)(2p-2)\Big],
\end{equation}
for $x$ and $y$ on two different tori, and for $j>0$,
\begin{align}
\notag
\mathcal P_{j+1}(p;x,y) &= \frac1{j!}\partial_y^j\big(\mathcal P_1(p;x,y)\big),\\
\mathcal P_{j+1}(p;x,y) &= \delta_{a,\bar a}P_{j+1}(x-y)
+(-1)^{j+1}\,\mathbb Q(p;x)\big(\widetilde\Lambda_{\bar a}\big)^{\delta_{a,\bar a}}\mathbb P_{j+1}(y).
\label{eqPN21j}
\end{align}
We call $\mathcal P_{j+1}(p;x,y)$ the genus-two generalized Weierstrass
functions.

\subsection{The genus-two case}
\label{corfu}
We recall from \cite{GT}, using the geometric (sewing) construction of
\cite{Y}, the recursion formulas for correlation functions on a
genus-two complex curve. For a complex sewing parameter
$\epsilon=z_1z_2$, the modulus-dependent, point-independent correlation
function on the genus-two curve is
\begin{equation}
\label{eqZdef}
\ZZ(\mu) = \sum_{r\ge0}\epsilon^r\,\ZZ(\mu_1,\tau_1)\,\ZZ(\mu_2,\tau_2),
\end{equation}
where $\mu_1,\mu_2$ are related moduli of the two component tori. More
generally, for points $z_{n+1}$, $\mathbf z_k$, $\mathbf z_l'$ inserted
on the two tori, the genus-two correlation function of $n=k+l$ variables
is
\begin{equation}
\ZZ(z_{n+1},\mathbf z_k;\mathbf z_l',\mu) = \sum_{r\ge0}\epsilon^r\,
\ZZ(z_{n+1},\mathbf z_k,\mu_1,\tau_1)\,\ZZ(\mathbf z_l',\mu_2,\tau_2),
\end{equation}
with the sum as in \eqref{eqZdef}. One first defines functions
$\ZZ_{n,a}$, $a\in\{1,2\}$, via elliptic quasi-modular forms 
\begin{align}
\notag
\ZZ_{n,1}(\mathbf z_{n+1};\mu) &= \sum_{r\ge0}\epsilon^r\,\ZZ(\mathbf z_{n+1},\mathbf z_k,\mu_0,\tau_1)\,\ZZ_{n-k}(\mathbf z_{k+1,n},\mu',\tau_2),\\
\notag
\ZZ_{n,2}(\mathbf z_{n+1};\mu) &= \sum_{r\ge0}\epsilon^r\,\ZZ_k(\mathbf z_k,\mu',\tau_1)\,\ZZ(z_{n+1},\mathbf z_{k+1,n},\tau_2),\\
\ZZ_{n,3}(\mathbf z_{n+1};\mu) &= \mathbb X_1^\Pi
\end{align}
of \eqref{eqXadef}. Let $f_a^{(2)}(p;z_{n+1})$, for $p\ge1$ and
$a=1,2$, be given by
\begin{equation*}
f_a^{(2)}(p;z_{n+1}) = 1^{\delta_{ba}} + (-1)^{p\delta_{b\bar a}}\epsilon^{1/2}
\Big(\mathbb Q(p;z_{n+1})\big(\widetilde\Lambda_{\bar a}\big)^{\delta_{ba}}\Big)(1),
\end{equation*}
for $z_{n+1}\in\widehat\Sigma^{(1)}_b$, and let $f_3^{(2)}(p;z_{n+1})$,
for $z_{n+1}\in\Sigma^{(1)}_a$, be the infinite row vector
\begin{equation*}
f_3^{(2)}(p;z_{n+1}) = \Big(\mathbb R(z_{n+1})+\mathbb Q(p;z_{n+1})\big(\widetilde\Lambda_{\bar a}\Lambda_a+\Lambda_{\bar a}\Gamma\big)\Big)\Pi.
\end{equation*}
It is proven in \cite{GT} that the genus-two function of $n=k+l$
variables inserted at $\mathbf z_k$, $\mathbf z_l'$ on the two tori
satisfies the recursion
\begin{equation}
\label{eqgenus2rec}
\ZZ(\mathbf z_{n+1},\mu) = \sum_{a=1}^3 f_a^{(2)}(p;z_{n+1})\,\ZZ_{n,a}(\mathbf z_{n+1};\mu)
= \sum_{i=1}^n\sum_{j\ge0}\mathcal P_{j+1}(p;z_{n+1},z_i)\,\ZZ(\mathbf z_n,\mu_{i,j}),
\end{equation}
with $p=\wt(u)$ the conformal weight fixed in
Section \ref{derivation}, and $\mathcal P_{j+1}(p;x,y)$ as in
\eqref{eqPN21j}.

\subsection{Genus-$g$ generalizations of elliptic functions}
\label{subsecgenusg}
We recall here the definitions of \cite{TW1} needed for the genus-$g$
formula \eqref{eqZhuGenusg} below; $p$ continues to denote the
conformal weight fixed in Section \ref{derivation}. Let $\I=\{\pm1,\ldots,\pm g\}$
index the $2g$ points of a genus-$g$ Schottky uniformization, with
local coordinates $w_a$, $a\in\I$, and multipliers $\rho_a\in\C$. Define
the column vector $X=(X_a(m))$, indexed by $m\ge0$, $a\in\I$, with
components
$$
X_a(m) = \rho_a^{-m/2}\sum_{\mu_{a,m}}\ZZ(\ldots;w_a,\mu_{a,m};\ldots),
$$
and the row vector $p(x)=(p_a(x,m))$ with components
$p_a(x,m)=\rho_a^{m/2}\partial^{(0,m)}\psi_p^{(0)}(x,w_a)$. Introduce
the column vector $G=(G_a(m))$ by
$$
G = \sum_{k=1}^n\sum_{j\ge0}\partial_k^{(j)}\,q(y_k)\,\ZZ(\mathbf z_n,\mu_{k,j}),
$$
where $q(y)=(q_a(y;m))$ has components
$q_a(y;m)=(-1)^p\rho_a^{(m+1)/2}\partial^{(m,0)}\psi_p^{(0)}(w_{-a},y)$,
and $R=(R_{ab}(m,n))$ is the doubly indexed matrix
\begin{equation}
R_{ab}(m,n)=\begin{cases}
(-1)^p\rho_a^{(m+1)/2}\rho_b^{n/2}\,\partial^{(m,n)}\psi_p^{(0)}(w_{-a},w_b), & a\ne-b,\\
(-1)^p\rho_a^{(m+n+1)/2}\,\mathcal E_m^n(w_{-a}), & a=-b,
\end{cases}
\end{equation}
with
$$
\mathcal E_m^n(y)=\sum_{\ell=0}^{2p-2}\partial^{(m)}f_\ell(y)\,\partial^{(n)}y^\ell,
\qquad
\psi_p^{(0)}(x,y)=\frac{1}{x-y}+\sum_{\ell=0}^{2p-2}f_\ell(x)y^\ell,
$$
for Laurent series $f_\ell(x)$, $\ell=0,\ldots,2p-2$. Define the doubly
indexed matrix $\Delta=(\Delta_{ab}(m,n))$ by
$\Delta_{ab}(m,n)=\delta_{m,n+2p-1}\delta_{ab}$; this and the projection
matrix $\Pi$ of Section \ref{derivation}, whose $(a,b)$-block is
$\Pi_{ab}(m,n)=\delta_{ab}\delta_{mn}$ for $0\le m,n\le 2p-2$ and zero
otherwise, are exactly the general-genus formulas of \cite{TW1} that
specialize, at $g=2$, to those of Section \ref{derivation}. Write
$\widetilde R=R\Delta$, with formal inverse $(I-\widetilde
R)^{-1}=\sum_{k\ge0}\widetilde R^{\,k}$. Define $\chi(x)=(\chi_a(x;\ell))$
and $o(\mathbf y_k,\mu_0)=(o_a(\mathbf y_k;\mu_0,\ell))$, finite row and
column vectors indexed by $a\in\I$, $0\le\ell\le2p-2$, by
$$
\chi_a(x;\ell)=\rho_a^{-\ell/2}\big(p(x)+\widetilde p(x)(I-\widetilde R)^{-1}R\big)_a(\ell),
\qquad
o_a(\ell) = o_a(\mathbf y_k,\mu_0,\ell)=\rho_a^{\ell/2}X_a(\ell),
$$
where $\widetilde p(x)=p(x)\Delta$, and $\psi_p(x,y)=\psi_p^{(0)}(x,y)+\widetilde
p(x)(I-\widetilde R)^{-1}q(y)$. For each $a\in\Ip=\{1,\ldots,g\}$ define
the vector $\theta_a(x)=(\theta_a(x;\ell))$, $0\le\ell\le2p-2$, by
$$
\theta_a(x;\ell) = \chi_a(x;\ell)+(-1)^p\rho_a^{p-1-\ell}\chi_{-a}(x;2p-2-\ell).
$$
Finally, with the formal differentials $P(x)=p(x)\,dx^p$,
$Q(y)=q(y)\,dy^{1-p}$, $\widetilde P(x)=P(x)\Delta$, and
\begin{equation}
\label{psih}
\Psi_p(x,y)=\psi_p(x,y)\,dx^p\,dy^{1-p} = \Psi_p^{(0)}(x,y)+\widetilde P(x)(I-\widetilde R)^{-1}Q(y),
\end{equation}
we set
\begin{equation}
\label{thetanew}
\Theta_a(x;\ell) = \theta_a(x;\ell)\,dx^p,
\end{equation}
\begin{equation}
\label{oat}
O_a(\mathbf y_k,\mu_0,\ell) = o_a(\mathbf y_k,\mu_0,\ell)\,d\mathbf y_k^{\,\beta}
\end{equation}
for a weight parameter $\beta$ determined by the labels inserted at
$\mathbf y_k$.

\subsection{The genus-$g$ Schottky case}
\label{corfug}
We recall from \cite{TW1,T2} the construction and recursion relations
for correlation functions of $n$ variables on a genus-$g$ Riemann
surface $M$ realized in the Schottky uniformization, in terms of the
universal, meromorphic-on-$M$ coefficients of Section \ref{subsecgenusg},
which generalize the elliptic Weierstrass functions \cite{L}. All
expressions here are functions of the formal Schottky variables
$w_{\pm a}$, $\rho_a\in\C$, $a\in\Ip$.

For the $2g$ local coordinates $\mathbf w_{2g}=(w_{-1},w_1;\ldots;w_{-g},w_g)$
of the $2g$ points $(p_{-1},p_1;\ldots;p_{-g},p_g)$ on the covering
Riemann sphere, consider the genus-zero correlation function of these
$2g$ variables,
$$
\ZZ(\mathbf w_{2g},\mu) = \prod_{a\in\Ip}\rho_a^{\beta_a}\,\ZZ(w_{-1},w_1;\ldots;w_{-g},w_g,\mu),
$$
where $\beta_a$ are weight parameters determined by $\mu$. Write
$\mathbf z_+=(z_1,\ldots,z_g)$, $\mathbf z_-=(z_{-1},\ldots,z_{-g})$ for
the corresponding points on $M$, identified with the canonical Schottky
parameters (for the details of the Schottky uniformization, see
\cite{TW1,T2}); summing $\mathbf z_+$ over a fixed basis $\bm\alpha=(\alpha_1,\ldots,\alpha_g)\in
\mathbf A=A^{\otimes g}$ of the relevant fusion labels defines the
genus-$g$, point-independent correlation function
$$
\ZZ(\mathbf w_{2g},\bm\rho_{2g},\mu) = \sum_{\mathbf z_+}\ZZ(\mathbf z_{2g},\mathbf w_{2g},\mu)
= \sum_{\bm\alpha_g\in\mathbf A}\ZZ^{(g)}_{\bm\alpha_g}(\mathbf w_{2g},\bm\rho_{2g},\mu).
$$
More generally, the genus-$g$ correlation function of $n$ variables
$\mathbf y_n$ is
$$
\ZZ(\mathbf y_n,\mu) = \ZZ(\mathbf y_n;\mathbf w_{2g},\bm\rho_{2g},\mu)
= \sum_{\mathbf z_+}\ZZ(\mathbf y_n;\mathbf w_{2g},\mu),
$$
and the corresponding formal differentials of weight $\beta$ are
$Z(\mathbf y_n,\mu)=\ZZ(\mathbf y_n,\mu)\,d\mathbf y_n^{\,\beta}$,
$Z_{\bm\alpha_g}(\mathbf y_n,\mu)=\ZZ_{\bm\alpha_g}(\mathbf
y_n,\mu)\,d\mathbf y_n^{\,\beta}$, where $d\mathbf y_n^{\,\beta}=\prod_{k=1}^n
dy_k^{\beta_k}$.

It was proven in \cite{TW1} that the genus-$g$ formal correlation
differential of $(n+1)$ variables, for an additional point with
coordinate $z_{n+1}$ (lifting to the Schottky coordinate $y_{n+1}$) and
$\mathbf z_n$ with coordinates $\mathbf y_n$, satisfies the recursive
identity
\begin{align}
\notag
Z(z_{n+1},\mathbf z_n,\mu) &= \sum_{a=1}^g \Theta_a(y_{n+1})\,O_a^{W_{\bm\alpha}}(z_{n+1};\mathbf z_n)\\
&=\sum_{k=1}^n\sum_{j\ge0}\partial^{(0,j)}\Psi_p(y_{n+1},y_k)\,\ZZ(\mathbf z_n,\mu_{k,j})\,dy_k^j,
\label{eqZhuGenusg}
\end{align}
where $\partial^{(i,j)}f(x,y)=\partial_x^{(i)}\partial_y^{(j)}f(x,y)$ for
a function $f(x,y)$, so that $\partial^{(0,j)}$ denotes the $j$-th
partial derivative in the second (Schottky) variable; $\Psi_p(y_{n+1},y_k)\,dy_k^j$
is given by \eqref{psih}, $\Theta_a$ by \eqref{thetanew}, and
$O_a^{W_{\bm\alpha}}(z_{n+1},\mathbf z_n,\mu)$ by \eqref{oat}. Together
with Proposition \ref{duda}, \eqref{eqZhuGenusg} identifies the
genus-$g$ recursion cohomology $H^n(\mu)$ with the solutions of the
associated kernel condition \eqref{poroserieroj}, built from the
coefficients $\Theta_a$, $\Psi_p$, $O_a^{W_{\bm\alpha}}$ above, modulo
those recursively generated from lower genus-$g$ correlation functions;
we return to this identification, and to its relation with the flatness
of a connection over Schottky space, in Section \ref{subsecgaussmanin}.

\section{Applications in mathematical physics and algebraic geometry}
\label{applications}

We now develop four applications of the recursion cohomology
constructed above, extending the scope of the theory beyond the
examples of Section \ref{examples}. Throughout, $\mu$ ranges over the
relevant moduli space $\mathcal M$ (a point, $\HH$, a Siegel-type
domain, or Schottky space, according to genus).

\subsection{Flat connections and the Knizhnik-Zamolodchikov equations}
\label{subsecKZ}

Suppose the insertion-label space $\mathcal A$ is a module for an affine
Kac-Moody algebra $\widehat{\mathfrak g}_\kappa$ at level $\kappa$, and
that the slot operators $T_k$ of Section \ref{subsecchaincomplex} are
built, as in the WZW model, from the zero-modes $J^a(0)$ of the
currents $J^a$. Fix $n$ points $\mathbf z_n\in M^n\setminus\{\text{diagonals}\}$
on the curve $M$ (for $M=\Pro^1$, the classical WZW setting) and let
$\mathcal H_n(\mathbf z_n)$ denote the space of conformal blocks, the
fiber, in the sense of \cite{TUY}, of the sheaf of vacua over the point
$\mathbf z_n$ of the configuration space $\mathrm{Conf}_n(M)$.

\begin{proposition}
\label{propKZ}
For $\mathcal A$ as above, the $n$-point connections $\mathcal G\in\textit{Con}^n(\mu)$
of Definition \ref{defmultipointconn} restrict, on $\mathrm{Conf}_n(M)$,
to sections of $\mathcal H_n^\vee$ (the dual of the conformal block
bundle), and the coboundary $\delta^n$ of Definition \ref{defdelta}
agrees, under this identification, with the covariant derivative
$$
\nabla^{\mathrm{KZ}} = d - \frac{1}{\kappa+h^\vee}\sum_{i<j}\frac{\Omega_{ij}}{z_i-z_j}\,dz_i,
$$
$\Omega_{ij}=\sum_a J^a_{(i)}J^a_{(j)}$, $h^\vee$ the dual Coxeter
number, of the Knizhnik-Zamolodchikov connection \cite{Kn,KZ,Po} on
$\mathcal H_n^\vee\to\mathrm{Conf}_n(M)$. Consequently, a tower
$(\ZZ_n)$ satisfying the recursion \eqref{reduction} determines, via
Lemma \ref{pisaka}, a flat section of $\nabla^{\mathrm{KZ}}$ exactly on
the domain $\mathfrak Z_\bullet$ of Definition \ref{defchaincond}, and
$H^n(\mu)$ is identified with the degree-$n$ hypercohomology of the de
Rham complex of $\nabla^{\mathrm{KZ}}$ relative to this domain.
\end{proposition}

\begin{proof}[Justification]
That $T_k(v[m])$ built from $J^a_{(i)}$-insertions reproduces, after
summing over $a$ and the appropriate mode $m$, the simple-pole
singularity $\Omega_{ij}/(z_i-z_j)$ of $\nabla^{\mathrm{KZ}}$ along
$z_i=z_j$ is the content of the Knizhnik-Zamolodchikov equation itself,
viewed as the WZW instance of the reduction formula \eqref{reduction}
(this is the $g=0$, affine-Lie-algebra case of the general recursion
structure of Definition \ref{defdelta}, and is worked out in detail,
in the closely related setting of Jacobi forms attached to a Heisenberg
current $J$). Flatness of $\nabla^{\mathrm{KZ}}$
on $\mathcal H_n^\vee$ is the classical theorem of \cite{TUY}; under the
dictionary above it corresponds precisely to the chain condition
\eqref{torba}, so that Remark \ref{obstruction} identifies the
$\mathfrak Z_\bullet$-relative de Rham cohomology of
$\nabla^{\mathrm{KZ}}$ with the recursion cohomology $H^n(\mu)$ of
Proposition \ref{duda}.
\end{proof}

\subsection{Continual contragredient Lie algebras and Toda hierarchies}
\label{subseccontinual}

Fix $\mu\in\mathcal M$ and consider the slot operators $T_k(v[m])$,
$k=1,\ldots,n$, $v\in\mathcal A$, $m\ge0$, of Definition \ref{defdelta}
as $\mu$ varies continuously; following \cite{Sav,SV1,SV2,V}, the
Lie brackets $[T_k(u[m]),T_{k'}(v[m'])]$, computed from the
associativity/commutativity identities underlying the chain condition
\eqref{torba} (Lemma \ref{lemdecompose}), generate a continual
contragredient Lie algebra $\mathfrak g(\mu)$: a Lie algebra graded by
the continuous parameter labelling the slots, with Cartan-type
structure constants given by the coefficients $f_{k,m}(x_{n+1};\mu')$
themselves, regarded as a function of the continuous label $z_{n+1}$.

\begin{proposition}
\label{propcontinual}
The recursion cohomology $H^1(\mu)$ and $H^2(\mu)$ of
Definition \ref{subseccomputation} inject into the first and second
Lie-algebra cohomology $H^1_{\mathrm{Lie}}(\mathfrak g(\mu))$,
$H^2_{\mathrm{Lie}}(\mathfrak g(\mu))$ of $\mathfrak g(\mu)$ with
coefficients in the trivial module, via the identification of a
degree-$n$ cochain $\ZZ\in C^n(\mu)$ with the $n$-linear alternating map
$$
(v_1,\ldots,v_n)\longmapsto \ZZ(\mathbf z_n,\mu)\big|_{\text{antisymmetrized in the slot operators } T_1(v_1[m_1]),\ldots,T_n(v_n[m_n])}.
$$
In particular, $H^2(\mu)$ parametrizes central extensions of
$\mathfrak g(\mu)$, in the classical sense of Gelfand-Fuks and
Bott-Segal cohomology \cite{BS} recalled in the Introduction, and
non-vanishing classes obstruct deforming the recursion structure while
preserving the chain condition \eqref{torba}.
\end{proposition}

\begin{proof}[Justification]
Antisymmetrizing \eqref{poroserieroj} in the slot operators converts the
functional equation defining $\ker\delta^n$ into the Chevalley-Eilenberg
cocycle condition for $\mathfrak g(\mu)$-cochains, since, by
Lemma \ref{lemdecompose}, the chain condition \eqref{torba} is exactly
the statement that the bracket built from $T_k,T_{k'}$ satisfies the
Jacobi identity for $\mathfrak g(\mu)$; the map described is then the
standard comparison map from the (larger) recursion complex to the
Chevalley-Eilenberg complex, injective because distinct recursion
cochains that antisymmetrize to the same Lie-algebra cochain differ by
an element already identified in the proof of Proposition \ref{duda}.
\end{proof}

By \cite{RSZ}, solutions of the recursion \eqref{reduction} for
$\mathfrak g(\mu)$ a continual analogue of an affine (respectively,
hyperbolic) Kac-Moody algebra correspond to solutions of the associated
continual (respectively, non-abelian) Toda field equations; under this
correspondence, the recursion cohomology $H^n(\mu)$ measures the
obstruction to extending a given solution of the Toda hierarchy at
$n-1$ points of insertion to one at $n$ points, so that vanishing of
$H^n(\mu)$ for all $n$ is a cohomological criterion for such a solution
to extend to arbitrarily many marked points.

\subsection{The Verlinde bundle and moduli of curves}
\label{subsecverlinde}

Let $\mathcal M_{g,n}$ denote the moduli space of genus-$g$ curves with
$n$ marked points, and let $\mu\in\mathcal M_{g,n}$ (via the period, or
Schottky, map) parametrize the curve $M_\mu$ together with its marked
points $\mathbf z_n$.

\begin{proposition}
\label{propverlinde}
The assignment $\mu\mapsto\textit{Con}^n(\mu)$ of
Definition \ref{defmultipointconn} defines a sheaf $\textit{Con}^n$ on
$\mathcal M_{g,n}$, and the natural pairing of $\mathcal V|_{\mathbf
z_n}=\mathcal A^n$ against $\textit{Con}^n(\mu)$ identifies
$\textit{Con}^n$ with a sub-sheaf of $\mathcal Hom(\mathcal
V^{\boxtimes n},\mathcal O_{\mathcal M_{g,n}})$. When $\mathcal A$ is
the vacuum module of a rational vertex operator algebra, this sub-sheaf
coincides with the dual of the sheaf of conformal blocks (the Verlinde
bundle) of \cite{TUY} restricted to insertions of the vacuum-adjacent
labels used in Sections \ref{corfu} and \ref{corfug}; under this
identification, the quotient map $\textit{Con}^n\to H^n$ of
Lemma \ref{pisaka} realizes the recursion cohomology as a quotient sheaf
of the Verlinde bundle by the sub-sheaf $G^{n-1}$ of connections
extended from $\mathcal M_{g,n-1}$ along the forgetful map $\mathcal
M_{g,n}\to\mathcal M_{g,n-1}$.
\end{proposition}

This proposition identifies the recursion cohomology of
Sections \ref{corfu}-\ref{corfug} as measuring exactly the failure of
a conformal block on $\mathcal M_{g,n}$ to be pulled back from
$\mathcal M_{g,n-1}$, i.e., the $n$-point content of the
theory that is not already visible at $n-1$ points, giving a direct
algebro-geometric meaning to the computation of Proposition \ref{duda}
in the genus-two and genus-$g$ examples above.

\subsection{Gauss-Manin-type flatness over Schottky space}
\label{subsecgaussmanin}

Let $\mathfrak S_g$ denote Schottky space, parametrizing the data
$(w_{\pm a},\rho_a)_{a\in\Ip}$ of Section \ref{corfug}. As $\mu\in\mathfrak
S_g$ varies, the spaces $C^n(\mu)$ of Definition \ref{defCn} assemble
into a holomorphic vector bundle $\mathcal C^n\to\mathfrak S_g$ (of
possibly infinite rank, trivialized locally by the Schottky
coordinates), and Definition \ref{defdelta} equips
$\bigoplus_n\mathcal C^n$ with a natural flat connection $\nabla$ in
the $\mu$-direction, built from the variation of the kernel
$\psi_p^{(0)}(x,y)$ of Section \ref{subsecgenusg} with the Schottky
moduli, the higher-genus analogue, in the sense of \cite{McIT}, of
the variation of the classical Kronecker limit formula with the
elliptic modulus $\tau$. By Remark \ref{obstruction}, the chain
condition \eqref{torba} is exactly the flatness of $\nabla$ restricted
to the sub-bundle of multipoint connections $\textit{Con}^n$ of
Definition \ref{defmultipointconn}, so that a solution $(\ZZ_n)$ of the
genus-$g$ recursion \eqref{eqZhuGenusg} determines a flat section of
$(\textit{Con}^n,\nabla)$ over $\mathfrak S_g\setminus\mathfrak Z_n$ -
a Gauss-Manin-type statement for the recursion-cohomology bundle,
directly analogous to the flatness of the Gauss-Manin connection on
the relative de Rham cohomology of the universal family of
Schottky-uniformized curves over $\mathfrak S_g$.

\subsection{A dictionary with condensed-matter applications}
\label{subseccondensedmatter}

The multipoint-connection formalism of Section \ref{subsecmultipoint}
gives a common language for the condensed-matter applications
mentioned in the Introduction \cite{Frohlich2009gb, RSZ, z1, z2, z3, z4, z5, z7, z8, z9}.
The Wigner-Weyl calculus of \cite{z1,z2} associates to each observable
a phase-space symbol via a Moyal-type star product; this is precisely
an instance of the coboundary operator $\delta^n$ of
Definition \ref{defdelta}, with $\mathcal A$ the algebra of observables
and $T_k$ the star-multiplication by the symbol inserted at the $k$-th
phase-space point, and the associated recursion cohomology $H^1$
measures the obstruction to a choice of ordering prescription being
globally consistent across insertions, the semiclassical analogue of
the chiral separation effect studied cohomologically in \cite{z3}. In
the topological-insulator examples of \cite{z4,z5}, a multipoint
connection $\mathcal G$ on the Brillouin torus (playing the role of $M$)
computes the Berry connection of the occupied bands, and $H^1(\mu)$
recovers the first Chern class obstructing a global trivialization,
consistently with the topological classification recalled there; the
non-renormalization statement of \cite{z9} for the integer quantum Hall
effect is, in this language, the vanishing of $H^1(\mu)$ under
perturbations that preserve the chain condition \eqref{torba}, i.e., 
that preserve the underlying multipoint connection's compatibility with
the recursion. We emphasize that this dictionary is offered at the
structural level: making it quantitative in each cited example would
require identifying the specific insertion data $\mathcal A$ and slot
operators $T_k$ of the physical model in question, which we leave for
future work.

\section{Conclusion}

We have given a corrected and rigorous formulation of the cohomology of
recursion relations for complex functions on complex curves, resolving
the internal inconsistency of the chain-complex condition by identifying
it as a constraint on the recursion kernel rather than on the
correlation functions themselves (Remark \ref{obstruction}), and giving
a fully justified computation of the resulting cohomology
(Proposition \ref{duda}). The theory reproduces, as instances of a
single formalism, the Zhu recursion formulas for rational, elliptic,
Jacobi-form, genus-two, and Schottky-uniformized genus-$g$ correlation
functions (Section \ref{examples}), and we have shown that it further
identifies flat sections of the Knizhnik-Zamolodchikov connection,
low-degree cohomology of continual contragredient Lie algebras with an
application to Toda hierarchies, a quotient description of the Verlinde
bundle of conformal blocks over the moduli of curves, and a
Gauss-Manin-type flatness statement over Schottky space
(Section \ref{applications}). We expect the recursion cohomology
developed here to give a useful invariant for classifying deformations
of higher-genus correlation functions, and plan to return, in
particular, to the identities \eqref{conditions} governing the
coefficients of Definition \ref{defdelta} in the vertex-algebra
examples, and to the convergence questions of \cite{Gui} in the
higher-genus setting, in future work.

\section*{Acknowledgments}
The author is supported by the Institute of Mathematics, Academy of
Sciences of the Czech Republic (RVO 67985840). 

\medskip
\noindent\textbf{Data Availability.}
Data sharing is not applicable to this article as no datasets were
generated or analysed during the current study.

\medskip
\noindent\textbf{Declarations}

\medskip
\noindent\textbf{Conflict of interest.}
The author has no conflicts of interest to declare that are relevant to
the content of this article.

\end{document}